\documentclass[12pt]{amsart}
\usepackage[centertags]{amsmath}
\usepackage{latexsym}
\usepackage{amssymb}

\DeclareMathOperator{\ad}{ad}

\DeclareMathOperator{\End}{\rm End}

\DeclareMathOperator{\bP}{\boldsymbol{\mathsf{P}}}
\DeclareMathOperator{\bS}{\boldsymbol{\mathsf{S}}}
\DeclareMathOperator{\bT}{\boldsymbol{\mathsf{T}}}
\DeclareMathOperator{\bTS}{\boldsymbol{\mathsf{T}}\boldsymbol{\mathsf{S}}}

\newcommand{\N}{\mathbb N}
\newcommand{\Q}{\mathbb Q}
\newcommand{\Z}{\mathbb Z}

\newcommand{\I}{\mathrm{i}}

\newtheorem{dummy}{Dummy}

\numberwithin{equation}{section}

\newtheorem{lemma}[dummy]{Lemma}

\newtheorem{theorem}[dummy]{Theorem}

\theoremstyle{definition}

\newtheorem{example}[dummy]{Example}

\theoremstyle{remark}
\newtheorem{rem}[dummy]{Remark}

\begin{document}

\bibliographystyle{amsalpha}
\author{Sandro Mattarei}
\email{mattarei@science.unitn.it}
\address{Dipartimento di Matematica\\
  Universit\`a degli Studi di Trento\\
  via Sommarive 14\\
  I-38100 Povo (Trento)\\
  Italy}
\title{Engel conditions and symmetric tensors}
\begin{abstract}
In a recent study of Engel Lie rings, Serena Cical\`{o} and Willem de Graaf
have given a practical set of conditions for an additively finitely generated Lie ring $L$ to
satisfy an Engel condition.
We present a simpler and more direct proof of this fact.
Then we generalize it to a
result in the language of tensor algebra, which can be applied to other contexts.
\end{abstract}
\thanks{The author acknowledges partial support from MIUR-Italy via PRIN
``Lie rings and algebras, groups, cryptography''.}
\maketitle

\section{Introduction}

Serena Cical\`{o} and Willem de Graaf~\cite{CicGra:I,CicGra:II} have recently developed algorithmic tools to investigate
finitely presented Lie rings.
They also have shown their effectiveness by applying them
to a computational study of some
finitely generated Lie rings satisfying an Engel condition.
One of the problems addressed there
is that of efficiently expressing an Engel condition
in terms of a generating set of the Lie ring as an additive group.
This turns out to be a substantial complication with respect to the more traditional case of Lie algebras over fields
of characteristic zero (or of sufficiently large characteristic).

Recall that a Lie ring $L$ is said {\em to satisfy the $n$-Engel condition}
(or {\em to be $n$-Engel} for short) if it satisfies the identity
$[\underbrace{x\cdots x}_{n}y]=0$ (using the right-normed convention for long Lie brackets, hence $[abc]:=[a[bc]]$ and so on).
This means that equality holds after substituting arbitrary elements of $L$ for $x$ and $y$.
If we assume, in addition, that $L$ is finitely generated as a Lie ring,
a celebrated result of Zelmanov~\cite{Zel:Burnside_odd} guarantees that L is nilpotent.
It follows that the additive group of $L$ is finitely generated,
say by elements $x_1,\ldots,x_m$.
In computational applications one faces the problem of efficiently expressing the condition that $L$ is $n$-Engel in terms
of linear combination of iterated Lie brackets in the additive generators $x_1,\ldots,x_m$.

If we worked instead with a Lie algebra $L$ over a field of characteristic larger than $n$,
a standard linearization (or polarization) trick would convert the $n$-Engel condition
into an equivalent multilinear identity, which could then be simply imposed on a basis for $L$ in all possible ways.
Explicitly, it is well known that a finite-dimensional such Lie algebra $L$ with basis $x_1,\ldots,x_m$ is $n$-Engel if and only if
\begin{equation}\label{eq:Engel_linearized}
\sum_{\sigma\in S_n}[x_{j_{\sigma(1)}}\cdots x_{j_{\sigma(n)}}y]=0,
\end{equation}
where $S_n$ denotes the symmetric group,
for all $1\le j_1\le\cdots\le j_n\le m$, and for each $y\in L$.

If $L$ is a Lie algebra over a field or ring where $n!$ is not invertible,
for example if $L$ is just a Lie ring (that is, a Lie algebra over $\Z$),
the multilinear conditions given in Equation~\eqref{eq:Engel_linearized} are still consequences
of the $n$-Engel condition, but they are generally not equivalent to it.
The following result, quoted from~\cite{CicGra:II}, gives a larger
set of conditions which is equivalent to the $n$-Engel condition in this setting,
with notation that we explain below.
The conditions expressed by Equation~\eqref{eq:Engel_linearized} are included in Equation~\eqref{eq:CG} as the case $s=n$.

\begin{theorem}[Theorem~18 in~\cite{CicGra:II}]\label{thm:CG}
Let $L$ be a Lie ring additively generated by $x_1,\ldots,x_m$.
Then $L$ satisfies the $n$-Engel condition
if and only if
\begin{equation}\label{eq:CG}
\sum_{\substack{k_1,\ldots,k_s>0\\ k_1+\cdots+k_s=n}}
[(x_{j_1}^{(k_1)}\cdots x_{j_s}^{(k_s)})^\ast y]=0
\end{equation}
holds for all $y\in L$, for
$1\le j_1\le\cdots\le j_s\le m$ and $1\le s\le n$.
\end{theorem}

Here the summand $[(x_{j_1}^{(k_1)}\cdots x_{j_s}^{(k_s)})^\ast y]$
can be defined as the result of substituting $x_{j_r}$ for $X_r$ into
the sum of all right-normed long Lie brackets, of weight $n+1$ and with $y$ as the right-most entry,
which are obtained by filling in $k_r$ copies of the symbol $X_r$, for each $r=1,\ldots,s$, in all possible ways.
(Note that the symbols $X_r$ are assumed distinct, but the elements $x_{j_r}$ which replace them need not be.)
In particular,  $[(x_{j_1}^{(k_1)}\cdots x_{j_s}^{(k_s)})^\ast y]$ is the sum of a number of long Lie brackets
equal to the multinomial coefficient
$\binom{n}{k_1,\ldots,k_s}=n!/(k_1!\cdots k_s!)$.
Because the Engel condition is linear in $y$, it is actually sufficient that Equation~\eqref{eq:CG}
holds for $y$ ranging over a set of additive generators for $L$.

Theorem~\ref{thm:CG} was deduced in~\cite{CicGra:II} from a previous coarser result of the authors~\cite[Theorem~14]{CicGra:I},
which instead of Equation~\eqref{eq:CG} had the much larger set of conditions
\begin{equation}\label{eq:CG_pm}
\sum_{\substack{k_1,\ldots,k_s>0\\ k_1+\cdots+k_s=n}}
p_{j_1}^{k_1}\cdots p_{j_s}^{k_s}
[(x_{j_1}^{(k_1)}\cdots x_{j_s}^{(k_s)})^\ast y]=0,
\end{equation}
to be satisfied for all choices of $p_{j_1},\ldots, p_{j_s}\in\{\pm1\}$.
We show in Section~\ref{sec:proof} how a simple observation entailing that the minus signs are superfluous leads to a direct and
much simpler proof of Theorem~\ref{thm:CG}.

Michael Vaughan-Lee has pointed out that our reduction
is analogous to one for groups, which is used in the GAP
implementation~\cite{GAP} of Werner Nickel's nilpotent quotient algorithm~\cite{Nickel:Engel}
to simplify testing the $n$-Engel condition in finitely generated nilpotent groups.
To the best of our knowledge, however, our observation has not been recorded in the literature.

Because the $n$-Engel condition can be more succinctly expressed as $(\ad x)^n=0$ for all $x\in L$,
it is apparent that Theorem~\ref{thm:CG} should be an instance of a general fact
which may occur in different contexts.
After introducing some symmetrization formulas in Section~\ref{sec:symmetrizations},
we formulate and prove such a result in Section~\ref{sec:general}, in the language of tensor algebra.
Our Theorem~\ref{thm:generators_P} gives a mildly redundant (see Section~\ref{sec:final}) set of generators for $\bP_n(M)$,
the $\Z$-module generated by the $n$-th tensor powers $\underbrace{x\otimes\cdots\otimes x}_n$
of all elements $x$ of $M$ in the tensor algebra $\bT(M)$, where $M$ is a finitely generated $\Z$-module.
There we also show how that implies Theorem~\ref{thm:CG}.

\section*{Acknowledgment}

The author is grateful to Michael Vaughan-Lee for useful discussions on an earlier version of this paper.

\section{A proof of Theorem~1}\label{sec:proof}

Before giving our proof of Theorem~\ref{thm:CG} we briefly digress on its original proof in~\cite{CicGra:I,CicGra:II}.
The longest part of the proof, given in~\cite{CicGra:II},
consists in showing that all the conditions in~\eqref{eq:CG_pm} are consequences of
their small subset in~\eqref{eq:CG}.
This is achieved by means of combinatorial arguments.
By contrast, the proof in~\cite[Theorem~14]{CicGra:I}
that the set of conditions~\eqref{eq:CG_pm} is equivalent to $L$ being $n$-Engel
consists of a straightforward calculation.

The presence of the additional signs $p_{j_1},\ldots,p_{j_s}$ in Equation~\eqref{eq:CG_pm}
arises in the proof by expressing an arbitrary element of $L$ in terms of sums and differences
of the given additive generators $x_1,\ldots x_m$, possibly with repetitions.
The core of our simplification is the observation that taking differences is superfluous,
and it suffices to test the $n$-Engel condition
just on sums of elements from $x_1,\ldots x_m$ (possibly with repetitions).

Note that Theorem~\ref{thm:CG} gives $\binom{m+s-1}{s}$ linear relations to be checked, for each possible value of $s$
and for each $y\in M$ (or rather in a set of additive generators for $L$).
Therefore, verifying the $n$-Engel condition for $L$ according to Theorem~\ref{thm:CG}
requires one to check a total of
\begin{equation}\label{eq:count}
\sum_{s=1}^n\binom{m+s-1}{s}
=\sum_{s=1}^n\left(\binom{m+s}{s}-\binom{m+s-1}{s-1}\right)
=\binom{m+n}{n}-1
\end{equation}
linear relations for each $y$.
These are much fewer than those given in~\eqref{eq:CG_pm}, which are
$\sum_{s=1}^n 2^s\binom{m+s-1}{s}\ge 2^n\binom{m+n-1}{n}=\bigl(2^nm/(m+n)\bigr)\binom{m+n}{n}$
for each $y$,
or half that number after the natural normalization $p_{j_1}=1$.
Reducing the number and the complexity of necessary verifications as much as possible is of utmost importance
in computational applications.
The linear relations required to verify the $n$-Engel condition according to Theorem~\ref{thm:CG}
are still redundant, but only slightly in number, as will be clear after Section~\ref{sec:general}.
We discuss this point further in Section~\ref{sec:final}.

Our simplification of the proof of Theorem~\ref{thm:CG} starts with the following simple observation
on finitely generated $\Z$-modules.

\begin{lemma}\label{lemma:tV}
Let $V$ be a finitely generated $\Z$-module, and $W$ a submodule.
If $W\subseteq tV$ for all positive integers $t$, then $W=0$.
\end{lemma}

\begin{proof}
Because $V$ is a direct sum of cyclic submodules, say $V=\bigoplus_i\Z v_i$,
and because $\bigcap_{t>0}t\Z=0$,
we have $\bigcap_{t>0}tV=0$, whence the conclusion.
\end{proof}

Of course it would be sufficient to let $t$ range over the prime powers in Lemma~\ref{lemma:tV},
or even on any set of integers which contains
multiples of arbitrarily high powers of each prime, but we will have no need for these variations.

\begin{lemma}\label{lemma:positive}
Let $L$ be a Lie ring whose additive group is generated by finitely many elements $x_1,\ldots,x_m$.
Then $L$ satisfies the $n$-Engel condition if and only if
$[\underbrace{x\cdots x}_{n}y]=0$ holds for all $y\in L$, and for all $x\in L$
which are nonnegative linear combinations of $x_1,\ldots,x_m$.
\end{lemma}

By {\em nonnegative linear combinations} we mean linear combinations with coefficients in the nonnegative integers.

\begin{proof}
We prove the sufficiency of the condition, its necessity being obvious.
In view of Lemma~\ref{lemma:tV}, where $V=L$ and $W$ is the additive subgroup of $L$ generated by all Lie brackets
$[\underbrace{x\cdots x}_{n}y]=0$, for $x,y\in L$, it is enough to show that $W\subseteq tL$ for every positive integer $k$.
The latter condition is equivalent to the quotient Lie ring $L/tL$ being $n$-Engel.
Thus, $L$ is $n$-Engel if (and only if) its quotient $L/tL$ is $n$-Engel for every positive integer $t$.
However, because $L/tL$ is finite (or, more precisely, because it is a torsion group), each of its elements
can be expressed as a nonnegative linear combination of the images of $x_1,\ldots,x_m$.
The conclusion follows.
\end{proof}

After these considerations, Theorem~\ref{thm:CG} can be proved in a similar fashion as~\cite[Theorem~14]{CicGra:I} was,
just by omitting the signs $p_{j_r}$.
Unfortunately, that proof is marred with notational errors
(where indices $j_r$ erroneously appear both as free variables and as bound to summation signs),
which we amend here for the reader's convenience.
To simplify the notation for sums such as those in Equations~\eqref{eq:CG} and~\eqref{eq:CG_pm}, as in~\cite{CicGra:II} we adopt the shorthand
$\sum_n$
for the summation symbol
\[
\sum_{\substack{k_1,\ldots,k_s>0\\ k_1+\cdots+k_s=n}},
\]
where the dummy variables $k_1,\ldots,k_s$ explicitly
appear in the summand omitted here.

\begin{proof}[Proof of Theorem~\ref{thm:CG}]
According to Lemma~\ref{lemma:positive},
the Lie ring $L$ is $n$-Engel if and only if
$[\underbrace{x\cdots x}_{n}y]=0$ holds for all $y\in L$, and for all $x\in L$
which are nonnegative linear combinations of elements from $x_1,\ldots,x_m$.
Any such linear combination can be written as a sum $x_{j_1}+\cdots+x_{j_s}$
of elements taken from $x_1,\ldots,x_m$, with $1\le j_1\le\cdots\le j_s\le m$.

We expand by linearity
\begin{align}\notag
[(x_{j_1}+\cdots+x_{j_s})&\cdots(x_{j_1}+\cdots+x_{j_s})y]
=
\\
\begin{split}\label{eq:x}
&=
\sum_{\substack{1\le r_1\le s}}
[(x_{j_{r_1}}^{(n)})^\ast y]
\\
&\quad+
\sum_{1\le r_1<r_2\le s}
\sum_n
[(x_{j_{r_1}}^{(k_1)}x_{j_{r_2}}^{(k_2)})^\ast y]
\\
&\quad\;\;\vdots
\\
&\quad+
\sum_{1\le r_1<\cdots<r_{s-1}\le s}
\sum_n
[(x_{j_{r_1}}^{(k_1)}\cdots x_{j_{r_{s-1}}}^{(k_{s-1})})^\ast y]
\\
&\quad+
\sum_n
[(x_{j_1}^{(k_1)}\cdots x_{j_s}^{(k_s)})^\ast y].
\end{split}
\end{align}
We give a more transparent version of this decomposition in Equation~\eqref{eq:X} .
If $L$ satisfies condition~\eqref{eq:CG} of Theorem~\ref{thm:CG},
for all $y\in L$, for
$1\le j_1\le\cdots\le j_s\le m$ and $s\ge 1$ (where the condition is void for $s>n$)
then the right-hand side of Equation~\eqref{eq:x} vanishes, and hence $L$ is $n$-Engel according to Lemma~\ref{lemma:positive}.

For the converse, if $L$ is $n$-Engel then the left-hand side of Equation~\eqref{eq:x} vanishes.
Working by induction on $s$ one shows that the last sum
in the right-hand side,
$\sum_n
[(x_{j_1}^{(k_1)}\cdots x_{j_s}^{(k_s)})^\ast y]$,
vanishes as well,
and hence Equation~\eqref{eq:CG} holds.
\end{proof}

The use of induction on $s$ in the `if' part of the proof can be replaced by a more direct argument.
This is based on a formula (essentially Equation~\eqref{eq:X_inverse}) which expresses the sum
$\sum_n
[(x_{j_1}^{(k_1)}\cdots x_{j_s}^{(k_s)})^\ast y]$
of Equation~\eqref{eq:CG}
as a linear combination of `Engel' Lie brackets $[\underbrace{z\cdots z}_{n}y]=0$, for suitable elements $z$ of $L$.
We describe that in the next section, in a more abstract setting.

\section{Symmetrizations}\label{sec:symmetrizations}

We prepare the ground for a generalization of Theorem~\ref{thm:CG} in the next section
by establishing some fairly general symmetrization formulas, which extend the traditional polarization trick over fields of characteristic zero.

We work in the non-commutative polynomial ring over $\Z$ (that is, the free associative ring, or $\Z$-algebra) on countably many
non-commuting indeterminates $X_1,X_2,X_3,\ldots$, and fix a positive integer $n$.
For each monomial
$X_{i_1}\cdots X_{i_n}$ of degree $n$ we define its symmetrization
\begin{equation}
(X_{i_1}\cdots X_{i_n})^\ast:=\sum_{\sigma\in S_n}X_{\sigma(i_1)}\cdots X_{\sigma(i_n)}.
\end{equation}
If all indeterminates $X_{i_1},\ldots,X_{i_n}$ are distinct, then the symmetrization
is a sum of $n!$ pairwise distinct monomials.
If some indeterminate appears several times, say $k$ times, among $X_{i_1},\ldots,X_{i_n}$, then
any given monomial in the right-hand side appears with coefficient a multiple of $k!$.
More formally, letting the symmetric group $S_n$ act
on the set of $n$-tuples via
$\sigma(i_1,\ldots,i_n)=(i_{\sigma^{-1}(1)},\ldots,i_{\sigma^{-1}(n)})$,
each monomial appears in $(X_{i_1}\cdots X_{i_n})^\ast$ with multiplicity equal to the order of the
stabilizer in $S_n$ of the $n$-tuple $(i_1,\ldots,i_n)$.
It is then convenient to define modified symmetrizations, as follows.
For an ordered partition $(k_1,\ldots,k_s)$ of $n$
(that is, an $s$-tuple of positive integers with sum $n$) and indices $i_1,\ldots,i_s$ we set
\begin{equation}\label{eq:def}
(X_{i_1}^{(k_1)}\cdots X_{i_s}^{(k_s)})^\ast
:=
\frac{1}{k_1!\cdots k_s!}(\underbrace{X_{i_1}\cdots X_{i_1}}_{k_1}\underbrace{X_{i_2}\cdots X_{i_2}}_{k_2}
\cdots\underbrace{X_{i_s}\cdots X_{i_s}}_{k_s})^\ast.
\end{equation}
Note that we are not assuming that the indices $i_1,\ldots,i_s$ are pairwise distinct in this definition,
but if they are (and only if they are) then
$(X_{i_1}^{(k_1)}\cdots X_{i_s}^{(k_s)})^\ast$
is a sum of
$\binom{n}{k_1,\ldots,k_s}=n!/(k_1!\cdots k_s!)$
pairwise distinct monomials.
More precisely, it is the sum of all distinct monomials obtained from $X_{i_1}^{k_1}\cdots X_{i_s}^{k_s}$
by permuting the factors.
In particular, we have
$(X_{i_1}^{(1)}\cdots X_{i_n}^{(1)})^\ast
=(X_{i_1}\cdots X_{i_n})^\ast$,
and
$(X_{i_1}^{(n)})^\ast
=X_{i_1}^n$.

In formal analogy with Equation~\eqref{eq:x}, and carrying over
the shorthand notation $\sum_n$
from Section~\ref{sec:proof},
we have
\begin{equation}
\begin{split}\label{eq:X}
(X_1+\cdots+X_s)^n
&=
\sum_{\substack{1\le r\le s}}
(X_{r}^{(n)})^\ast
\\
&\quad+
\sum_{1\le r_1<r_2\le s}
\sum_n
(X_{r_1}^{(k_1)}X_{r_2}^{(k_2)})^\ast
\\
&\quad\;\;\vdots
\\
&\quad+
\sum_{1\le r_1<\cdots<r_{s-1}\le s}
\sum_n
(X_{r_1}^{(k_1)}\cdots X_{r_{s-1}}^{(k_{s-1})})^\ast
\\
&\quad+
\sum_n
(X_1^{(k_1)}\cdots X_s^{(k_s)})^\ast,
\end{split}
\end{equation}
which amounts to sorting the resulting monomials according to the total number of distinct indeterminates which they contain.
Now we show how Equation~\eqref{eq:X} can be `inverted'.

For any nonempty subset $I$ of $\N=\{1,2,3,\ldots\}$ let $X_I$ be the sum of all distinct
monomials of degree $n$ formed with the indeterminates $X_i$ with $i\in I$, each appearing at least once.
Hence $X_I=0$ if $|I|>n$, being the sum of an empty set of monomials.
We also conveniently set $X_{\varnothing}=0$.
If $I=\{i_1,\ldots,i_{|I|}\}$ with $i_1<\cdots<i_{|I|}$, we have
\[
X_I:=\sum_{\substack{k_1,\ldots,k_{|I|}>0\\ k_1+\cdots+k_{|I|}=n}}
(X_{i_1}^{(k_1)}\cdots X_{i_{|I|}}^{(k_{|I|})})^\ast
\]
Then we have
\[
\biggl(\sum_{i\in I} X_i \biggr)^n
=
\sum_{J\subseteq I} X_{J}.
\]
Note that the special case $I=\{1,\ldots,s\}$ of this decomposition is formally analogous to Equation~\eqref{eq:x}.
M\"obius inversion (on the poset $\mathcal{P}(\N)$ of subsets of $\N$,
see~\cite[p.~463]{Jac:BAI}) yields
\[
X_I=
\sum_{J\subseteq I}
(-1)^{|I\setminus J|}
\biggl(\sum_{j\in J} X_j \biggr)^n.
\]
In our case of interest $I=\{1,\ldots,s\}$ this formula reads
\begin{equation}\label{eq:X_inverse}
\sum_n
(X_1^{(k_1)}\cdots X_s^{(k_s)})^\ast
=
(-1)^s
\sum_{J\subseteq\{1,\ldots,s\}}
(-1)^{|J|}
\biggl(\sum_{j\in J} X_j \biggr)^n,
\end{equation}
which achieves the desired inversion of Equation~\eqref{eq:X}.
The special case $s=n$ of Equation~\eqref{eq:X_inverse}, where the left-hand side reads
$(X_1\cdots X_n)^\ast$,
is well known, see~\cite[Chapter~I, \S8.2, Proposition~2]{Bourbaki:Algebra_I}
or~\cite[Equation~(1.1)]{Kostrikin:Burnside}.
In fact, Equation~\eqref{eq:X_inverse}
is an explicit and characteristic-free version of~\cite[Equation~(1.2)]{Kostrikin:Burnside}.

\begin{example}
When $n=4$, Equation~\eqref{eq:X_inverse} reads
\[
(X_1^{(3)}X_2^{(1)})^\ast+(X_1^{(2)}X_2^{(2)})^\ast+(X_1^{(1)}X_2^{(3)})^\ast
=
(X_1+X_2)^4-X_1^4-X_2^4
\]
for $s=2$, and
\begin{multline*}
(X_1^{(2)}X_2^{(1)}X_3^{(1)})^\ast+(X_1^{(1)}X_2^{(2)}X_3^{(1)})^\ast+(X_1^{(1)}X_2^{(1)}X_3^{(2)})^\ast
=\\=
(X_1+X_2+X_3)^4
-(X_1+X_2)^4-(X_1+X_3)^4-(X_2+X_3)^4
+X_1^4+X_2^4+X_3^4
\end{multline*}
for $s=3$.
Note that in an ordinary commutative polynomial ring the left-hand sides of these equations would read
$4X_1^3X_2+6X_1^2X_2^2+4X_1X_2^3$
and
$12X_1^2X_2X_3+12X_1X_2^2X_3+12X_1X_2X_3^2$.
\end{example}

We extend the notation for symmetrizations to an arbitrary associative ring by letting
$(x_{i_1}^{(k_1)}\cdots x_{i_s}^{(k_s)})^\ast$
be the result of evaluating
$(X_{i_1}^{(k_1)}\cdots X_{i_s}^{(k_s)})^\ast$
on elements $x_i$ of the ring.
Equations~\eqref{eq:X} and~\eqref{eq:X_inverse} then hold with elements $x_i$
in place of the indeterminates $X_i$.

\section{A generalization in terms of tensor algebra}\label{sec:general}

Let $A$ be a commutative ring, and let $M$ be an $A$-module.
Let $\bT(M)=\bigoplus_{n=0}^{\infty}\bT^n(M)$
be the tensor algebra of $M$, and
let $\bS(M)=\bigoplus_{n=0}^{\infty}\bS^n(M)$
be the symmetric algebra of $M$.
Thus, $\bS(M)$ is the quotient of $\bT(M)$ by the ideal generated by the elements $x\otimes y-y\otimes x$,
for $x,y\in M$,
and we have a natural morphism of graded algebras of $\bT(M)$ onto $\bS(M)$.
Recall that the symmetric group $S_n$ acts on
$\bT^n(M)$ by linear extension of the action
\[
\sigma(x_1\otimes x_2\otimes\cdots\otimes x_n)
=
x_{\sigma^{-1}(1)}\otimes x_{\sigma^{-1}(2)}\otimes\cdots\otimes x_{\sigma^{-1}(n)}
\]
on the set of pure tensors, for $x_1,x_2,\ldots,x_n\in M$ and $\sigma\in S_n$.
Let $\bS'_n(M)$ (as denoted in~\cite[Chapter~III, \S6.3]{Bourbaki:Algebra_I},
rather than $\bTS^n(M)$ as in~\cite[Chapter~IV, \S5.3]{Bourbaki:Algebra_II})
consist of the symmetric tensors of degree $n$, that is,
of the elements of $\bT^n(M)$ fixed under the action of $S_n$.
Thus, $\bS'_n(M)$ is an $A$-submodule of $\bT^n(M)$.

It is well known that $\bS_n(M)$ and $\bS'_n(M)$ are naturally isomorphic when $n!$ (that is, $n!\cdot 1_A$)
is invertible in $A$, but not otherwise.
In general, the symmetrization map
$\boldsymbol{s}\colon z\mapsto\sum_{\sigma\in S_n}\sigma z$
maps $\bT^n(M)$ onto an $A$-submodule $\bS''_n(M)$ of $\bS'_n(M)$,
and factors through the natural epimorphism $\bT^n(M)\to\bS^n(M)$,
producing a morphism $\boldsymbol{\bar s}\colon\bS^n(M)\to\bS'_n(M)$ of $A$-modules.
Because $\boldsymbol{s}(z)=n!\,z$ for $z\in\bS'_n(M)$,
when $n!$ is invertible in $A$ the morphism $\boldsymbol{\bar s}$ is surjective, and so $\bS''_n(M)=\bS'_n(M)$.
One can prove that in this case $\boldsymbol{\bar s}$ is injective as well,
and we refer to~\cite[Chapter~III, \S6.3]{Bourbaki:Algebra_I} for a proof.
However, $\bS''_n(M)$ differs from $\bS'_n(M)$ when $n!$ is not invertible in $A$,
for example when $A=\Z$, a case which we examine more closely below.

For arbitrary $A$ and $M$,
the $A$-module $\bS^n(M)$ is generated by the set
$\{x_1\cdots x_n\mid x_i\in M\}$,
and hence $\bS''_n(M)$ is generated by its image
$\{\boldsymbol{\bar s}(x_1\cdots x_n)\mid x_i\in M\}$ under the symmetrization map.
We introduce an intermediate object $\bP_n(M)$ between
$\bS''_n(M)$ and $\bS'_n(M)$, which is the $A$-submodule of $\bT^n(M)$ generated by the set
$\{\underbrace{x\otimes\cdots\otimes x}_n\mid x\in M\}$
of all $n$-th powers (in $\bT^n(M)$) of elements of $M$.
While it is clear that $\bP_n(M)\subseteq\bS'_n(M)$,
the inclusion $\bS''_n(M)\subseteq\bP_n(M)$ follows from
evaluating the special case $s=n$ of Equation~\eqref{eq:X_inverse}
on arbitrary elements $x_1,\ldots,x_n$ of $M$,
and thus expressing the symmetrization $\boldsymbol{\bar s}(x_1\cdots x_n)$
as an $A$-linear combination of $n$-th powers.

Now suppose that $A=\Z$ and that $M$ is a free $\Z$-module of finite rank.
Then so is $\bT^n(M)$, and hence so are its submodules
$\bS''_n(M)\subseteq\bP_n(M)\subseteq\bS'_n(M)$.
In fact, the latter three $\Z$-modules all have the same rank $\binom{m+n-1}{n}$,
which is the dimension of
$\bS''_n(\Q\otimes_\Z M)=\bP_n(\Q\otimes_\Z M)=\bS'_n(\Q\otimes_\Z M)$.
To see this, note that the morphism
$\boldsymbol{\bar s}\colon\bS^n(M)\to\bS''_n(M)$,
surjective by definition, yields a morphism
$\Q\otimes_\Z\bS^n(M)\to\Q\otimes_\Z\bS''_n(M)$,
also surjective by right-exactness of the functor $\otimes_\Z$.
However, if we canonically identify
$\Q\otimes_\Z\bS^n(M)$ with $\bS^n(\Q\otimes_\Z M)$,
and similarly $\Q\otimes_\Z\bT^n(M)$ with $\bT^n(\Q\otimes_\Z M)$
(see~\cite[Chapter~III, \S6.4]{Bourbaki:Algebra_I})
we conclude that
$\Q\otimes_\Z\bS''_n(M)=\bS''_n(\Q\otimes_\Z M)$.
Thus, $\bS''_n(M)$ is a full sublattice in
$\bS''_n(\Q\otimes_\Z M)$,
and hence so are the $\Z$-modules
$\bP_n(M)$ and $\bS'_n(M)$, as they contain $\bS''_n(M)$.
Our assertion on their common rank follows from the fact that
$\bS''_n(\Q\otimes_\Z M)$, being isomorphic to $\bS_n(\Q\otimes_\Z M)$ because $1/n!\in\Q$, has dimension $\binom{m+n-1}{n}$.

At this point it is natural to ask for bases of  $\bS''_n(M)$, $\bP_n(M)$ and $\bS'_n(M)$,
or at least for reasonably small generating sets.
Bases for $\bS''_n(M)$ and $\bS'_n(M)$ are easy to produce,
and later we will obtain a mildly redundant generating set for $\bP_n(M)$.

\begin{theorem}\label{thm:bases}
Let $M$ be a free $\Z$-module of finite rank $m$, and let $x_1,\ldots,x_m$ be a free basis.
Then the following statements hold.
\begin{enumerate}
\item
A basis of $\bS'_n(M)$ is given by the set of elements of $M$ of the form
\begin{equation*}
(x_{j_1}^{(k_1)}\cdots x_{j_s}^{(k_s)})^\ast
\end{equation*}
where $1\le j_1<\cdots<j_s\le m$,
$k_1,\ldots,k_s>0$ with $k_1+\cdots+k_s=n$,
and $s=1,\ldots n$.
\item
A basis of $\bS''_n(M)$ is given by the set of elements of $M$ of the form
\begin{equation*}
k_1!\cdots k_s!\,(x_{j_1}^{(k_1)}\cdots x_{j_s}^{(k_s)})^\ast
\end{equation*}
where $1\le j_1<\cdots<j_s\le m$,
$k_1,\ldots,k_s>0$ with $k_1+\cdots+k_s=n$,
and $s=1,\ldots n$.
\end{enumerate}
\end{theorem}

\begin{proof}
To prove assertion~(1), note first that the given elements of $\bS'_n(M)\subseteq\bT^n(M)$ are $\Z$-independent.
Now write an arbitrary element $z$ of $\bS'_n(M)$
(that is, a symmetric tensor in $\bT^n(M)$) as a $\Z$-linear combination of pure tensors in the basis elements.
If one such pure tensor has $s$ distinct basis elements
$x_{j_1},\ldots,x_{j_s}$ as factors, appearing with multiplicities $k_1,\ldots,k_s$ in some order,
then all the pure tensors obtained from it by permuting the factors must appear with the same coefficient in $z$.
Because the sum of all those factors equals $(x_{j_1}^{(k_1)}\cdots x_{j_s}^{(k_s)})^\ast$,
the desired conclusion follows.

To prove assertion~(2) note that the monomials $x_{j_1}^{k_1}\cdots x_{j_s}^{k_s}$ of total degree $n$
form a basis of $\bS_n(M)$, and apply the symmetrization map $\boldsymbol{\bar s}\colon\bS_n(M)\to\bS'_n(M)$,
which is injective because so is its extension
$\boldsymbol{\bar s}\colon\bS_n(\Q\otimes_\Z M)\to\bS'_n(\Q\otimes_\Z M)$, making the identifications described earlier.
\end{proof}

\begin{rem}\label{rem:order}
According to Theorem~\ref{thm:bases}, the index of $\bS'_2(M)$ in $\bS''_2(M)$
is given by the formula
\[
|\bS'_n(M)/\bS''_n(M)|=
\prod_{s=1}^n
\prod_{\substack{k_1,\ldots,k_s>0\\ k_1+\cdots+k_s=n}}
(k_1!\cdots k_s!)^{\binom{m}{s}}.
\]
For example, we have
\begin{align*}
|\bS'_2(M)/\bS''_2(M)|&=2^m,
\\
|\bS'_3(M)/\bS''_3(M)|
&=(2\cdot 3)^m\cdot(2\cdot 2)^{\binom{m}{2}}
=2^{m^2}3^m,
\\
|\bS'_4(M)/\bS''_4(M)|
&=(2^3\cdot 3)^m\cdot(2^4\cdot 3^2)^{\binom{m}{2}}\cdot(2^3)^{\binom{m}{3}}
=2^{(m^2+m+4)m/2}\cdot 3^{m^2},
\\
|\bS'_5(M)/\bS''_5(M)|
&=(2^3\cdot 3\cdot 5)^m\cdot(2^{10}\cdot 3^4)^{\binom{m}{2}}\cdot(2^9\cdot 3^3)^{\binom{m}{3}}\cdot(2^4)^{\binom{m}{4}}
\\&=2^{(m^2+3m+14)m^2/6}\cdot 3^{(m^2+m+4)m/2}\cdot 5^m.
\end{align*}
In general, $|\bS'_n(M)/\bS''_n(M)|$ is divisible only by primes not exceeding $n$.
Furthermore, the highest power of a prime $p$ which divides $|\bS'_p(M)/\bS''_p(M)|$ is $p^m$.
\end{rem}

Both assertions of Theorem~\ref{thm:bases} certainly hold in greater generality than as stated, but we have restricted
our attention to $\Z$-modules for a direct comparison with our result on $\bP_n(M)$ below,
where that assumption is more crucial.
In order to provide a generating set for $\bP_n(M)$ which is not too redundant
we will need the following more general version of Lemma~\ref{lemma:tV}.

\begin{lemma}\label{lemma:tV_extended}
Let $V$ be a finitely generated $\Z$-module, and $U,W$ submodules.
If $U+tV=W+tV$ for all positive integers $t$, then $U=W$.
\end{lemma}

\begin{proof}
Because $U+W+tV=U+tV$ we have
$(U+W)/U\subseteq t(V/U)$, for all positive integers $t$.
Now Lemma~\ref{lemma:tV} applies with $V/U$ instead of $V$, and shows that
$(U+W)/U=0$.
Hence $W\subseteq U$, and the desired conclusion follows by symmetry.
\end{proof}

\begin{theorem}\label{thm:generators_P}
Let $M$ be a free $\Z$-module of finite rank $m$, and let $x_1,\ldots,x_m$ be a free basis.
Then $\bP_n(M)$ is generated by the set of elements of $M$ of the form
\begin{equation*}
\sum_n
(x_{j_1}^{(k_1)}\cdots x_{j_s}^{(k_s)})^\ast
\end{equation*}
for $1\le j_1\le\cdots\le j_s\le m$ and $1\le s\le n$.
\end{theorem}

\begin{proof}
By evaluating Equation~\eqref{eq:X_inverse} on the elements $x_{j_1},\ldots,x_{j_s}$
(viewed as elements of degree one in the tensor algebra $\bT(M)$)
one obtains an expression for
$(x_{j_1}^{(k_1)}\cdots x_{j_s}^{(k_s)})^\ast$
as a linear combination of $n$-th powers of elements of $M$, thus proving that it belongs to $\bP_n(M)$.

Before dealing with the converse, we need a reduction similar to Lemma~\ref{lemma:positive}.
Let $U$ be the $\Z$-submodule of $\bP_n(M)$ generated by
the elements $\underbrace{y\otimes\cdots\otimes y}_n$
where $y$ ranges over the {\em nonnegative} linear combinations of $x_1,\ldots,x_m$
(defined in Section~\ref{sec:proof}).
We prepare for an application of Lemma~\ref{lemma:tV_extended} to show that $U=\bP_n(M)$.
Let $t$ be a positive integer.
If $x$ is any element of $M$, hence a $\Z$-linear combination of $x_1,\ldots,x_m$,
then we can write $x=y+tz$, for some $y,z\in M$ with $y$ a nonnegative linear combination of $x_1,\ldots,x_m$.
By expanding the $n$-th power of $x$ we find that
\[
\underbrace{x\otimes\cdots\otimes x}_n
\equiv
\underbrace{y\otimes\cdots\otimes y}_n
\pmod{t\bT^n(M)},
\]
and hence
$\bP_n(M)+t\bT^n(M)=U+t\bT^n(M)$.
Because this holds for all positive integers $t$,
Lemma~\ref{lemma:tV_extended} with $V=\bT^n(M)$ implies that $U=\bP_n(M)$, as desired.

The rest of the proof runs similar to the proof of Theorem~\ref{thm:CG}.
We only need to show that any element $\underbrace{y\otimes\cdots\otimes y}_n$,
with $y$ a nonnegative linear combination of $x_1,\ldots,x_m$,
can be expressed as a $\Z$-linear combination of the elements given in Theorem~\ref{thm:generators_P}.
This is achieved by expressing $y$ as a sum $y=x_{j_1}+\cdots+x_{j_s}$,
with $1\le j_1\le\cdots\le j_s\le m$,
and by evaluating
Equation~\eqref{eq:X}
on the elements $x_{j_1},\ldots,x_{j_s}$.
\end{proof}

Because Theorem~\ref{thm:generators_P} only provides a set of generators for $\bP_n(M)$ and,
differently from Theorem~\ref{thm:bases}, makes no claim of independence,
it remains true for any finitely generated $\Z$-module $M$
having $x_1,\ldots,x_m$ a set of generators.
In fact, that follows from the present version applied to a free $\Z$-module $M'$ with basis $x'_1,\ldots,x'_m$,
to be mapped onto $M$ in the obvious way.

\begin{rem}\label{rem:p}
Theorem~\ref{thm:generators_P} is really of interest only for $n>2$, because $\bP_2(M)=\bS'_2(M)$.
In fact, the generating set for $\bP_2(M)$ given in Theorem~\ref{thm:generators_P} includes
the basis of $\bS'_2(M)$ given in Theorem~\ref{thm:bases}.
Incidentally, the fact that $\bP_2(M)=\bS'_2(M)$ admits the following generalization:
if $p$ is a prime then $p$ does not divide
$\bS'_p(M)/\bP_p(M)$.
This follows from the fact that
$\bP_p(M/pM)=\bS'_p(M/pM)$,
a proof of which can be found in~\cite[Chapter~III, \S6, Exercise~5]{Bourbaki:Algebra_I}.
\end{rem}

To illustrate the flexibility of Theorem~\ref{thm:generators_P} we show how to deduce Theorem~\ref{thm:CG} from it.

\begin{proof}[Another proof of Theorem~\ref{thm:CG}]
To avoid confusion with the notation of this section,
denote the given generators of the additive group of $L$
with $x_1',\ldots,x_m'$, rather than $x_1,\ldots,x_m$ as in Theorem~\ref{thm:CG}.
The image $\ad L$  of the adjoint representation $\ad:L\to\End_\Z(L)$ of $L$
is then generated by $\ad x_1',\ldots,\ad x_m'$ as a $\Z$-submodule of $\End_\Z(L)$.
Let $M$ be a free $\Z$-module with free basis $x_1,\ldots,x_m$.
The map of $M$ to $\End_\Z(L)$ sending $x_j$ to $x_j'$
extends uniquely to a morphism $\bT(M)\to\End_\Z(L)$ of $\Z$-modules,
and its image is the enveloping algebra of $\ad L$ (that is, the smallest (unital associative) subalgebra of $\End_\Z(L)$
which contains $\ad L$).
This morphism maps $\bP_n(M)$ onto the $\Z$-submodule of $\End_\Z(L)$ generated by the set
$\{(\ad x)^n\mid x\in\ L\}$.
Hence $L$ is $n$-Engel if and only if $\bP_n(M)$ is mapped to zero.
According to Theorem~\ref{thm:generators_P},
this occurs if and only if the condition expressed in Equation~\eqref{eq:CG} is satisfied.
\end{proof}

Similar arguments allow applications of Theorem~\ref{thm:generators_P} to other contexts, such as
testing an associative ring for a nil condition (of given index $n$).

\section{Final comments}\label{sec:final}

\subsection{Redundance in Theorem~\ref{thm:bases}}

According to a calculation done in Equation~\eqref{eq:count},
the generating set for $\bP_n(M)$ provided by Theorem~\ref{thm:bases} has cardinality
$\binom{m+n}{n}-1=\frac{m+n}{m}\binom{m+n-1}{n}-1$.
Hence when $m$ and $n$ are large this number is larger than the cardinality of a basis by roughly a factor $1+n/m$.
In computational applications, where the rank $m$ is typically large compared with $n$,
this ratio is close to $1$.
However, there may be reasons other than a negligible gain in cardinality
for obtaining a basis of $\bP_n(M)$ from the set of generators given in Theorem~\ref{thm:bases}.

For example, the natural way of producing a basis of a free $\Z$-module from a set of generators
is by reducing a certain matrix to Hermite normal form.
At the end of Section~4 in~\cite{CicGra:II}, its authors describe
a specific way of doing this reduction in their more specialized context,
where the generators of our Theorem~5
play the role of equations as in Theorem~1.
Inspection of their results reveals that most of the equations resulting from this reduction
are much shorter than the original ones (that is, they contain fewer summands).
In light of this remark, Willem de Graaf agrees that the equations being shorter is more likely to be responsible for
the observed reduction in running times of their algorithm
than being marginally fewer in number
(as suggested in the last sentence of Section~4 in~\cite{CicGra:II}).

Given a specific value of $n$, reduction to Hermite normal form as described in~\cite{CicGra:II}
applies in our more general context
to produce a basis of $\bP_n(M)$ from the set of generators given in our Theorem~\ref{thm:generators_P}.
However, it is difficult to see how to give a general analysis of this reduction,
and thus describe the resulting basis of $\bP_n(M)$ for generic $n$.

\subsection{Other issues on an $n$-Engel test for Lie rings}
We should mention that the problem of efficiently checking whether a finitely generated Lie ring or algebra is $n$-Engel
has subtler issues, which we have deliberately disregarded in this paper.
In particular, in Section~\ref{sec:proof} we have made no use of the anticommutativity or the Jacobi relation
in the Lie ring, focussing only on a multilinear algebra aspect of the Engel condition, in preparation of our generalization.
As enlightening discussions with Michael Vaughan-Lee have made clear,
taking the Lie ring structure into account
allows, in practice, for very substantial reduction in the number of conditions necessary
to prove that a {\em given} finitely generated Lie ring $L$ is $n$-Engel, with respect to those given in Theorem~\ref{thm:CG}.

We mention only the simplest source of such a reduction as an example.
According to Zelmanov's result, a candidate $L$ to be $n$-Engel must be nilpotent for a start.
Thus, the generators $x_1,\ldots,x_m$ of $L$ as $\Z$-module can be assigned weights,
according to the latest term of the lower central series of $L$ which they belong to.
Assuming that these weights have been calculated, the various summands in Equation~\eqref{eq:CG} acquire weights accordingly.
Instances of the equation where all of these weights exceed the nilpotency class of $L$
are trivially satisfied in $L$ and, therefore, need not be checked.
This observation can be refined in various ways, and its value clearly depends on how the given Lie ring $L$ is specified.

\subsection{Other base rings}

We mentioned earlier that working with $\Z$-modules rather than modules over other commutative rings
is quite crucial for Theorem~\ref{thm:generators_P}.
We briefly elaborate on that point now.
We note that Theorem~\ref{thm:generators_P} remains valid if we replace the ring $\Z$ with any subring $A$ of $\Q$.
In fact, when $A=\Q$ we have $\bP_n(M)=\bS'_n(M)$,
and Theorem~\ref{thm:bases} (that is, its version over $\Q$) gives a better answer than Theorem~\ref{thm:generators_P}.
Now note that Lemma~\ref{lemma:tV}, and hence the more general Lemma~\ref{lemma:tV_extended},
continue to hold after replacing $\Z$ with a {\em proper} subring $A$ of $\Q$:
if $V$ is a finitely generated $A$-module, and $U,W$ are submodules
with $U+tV=W+tV$ for all positive integers $t$, then $U=W$.
Lemma~\ref{lemma:tV_extended} actually holds in greater generality, but
our present assumption that $A$ is a subring of $\Q$ implies that
$A/tA$ is cyclic (a quotient ring of $\Z$), which makes the proof of Theorem~\ref{thm:generators_P} work.

Rings which may be of interest for applications include rings of algebraic integers.
However, the following example illustrates how Theorem~\ref{thm:generators_P} does not extend, in general, to rings $A$
with non-cyclic quotients $A/tA$.

\begin{example}
Let $A=\Z[\I]$ be the ring of Gaussian integers, and let
$M$ be a free $A$-module of rank $2$, with a basis $x_1,x_2$.
Theorem~\ref{thm:generators_P} gives a set of $9$ generators for
the $\Z$-module $\bP_3(\Z x_1+\Z x_2)$, which can be reduced to a basis
\begin{equation*}
(x_1^{(3)})^\ast,\quad
(x_2^{(3)})^\ast,\quad
(x_1^{(2)}x_2^{(1)})^\ast+(x_1^{(1)}x_2^{(2)})^\ast,\quad
2(x_1^{(1)}x_2^{(2)})^\ast.
\end{equation*}
These elements, however, generate a proper $A$-submodule of $\bP_3(M)$, to which the element
\begin{align*}
(1+\I)(x_1^{(2)}x_2^{(1)})^\ast
&=(x_1+x_2)\otimes(x_1+x_2)\otimes(x_1+x_2)
\\
&\quad+(x_1+\I x_2)\otimes(x_1+\I x_2)\otimes(x_1+\I x_2)
\\
&\quad-2 x_1\otimes x_1\otimes x_1
-(1-\I) x_2\otimes x_2\otimes x_2
\end{align*}
of $\bP_3(M)$ does not belong.
\end{example}

One possible way to deal with a ring $A$ which is finitely generated as a $\Z$-module, say by $r$ elements,
is forgetting the $A$-module structure of $M$ and viewing $M$ as a $\Z$-module on $rm$ generators.
Theorem~\ref{thm:generators_P} applies to the latter, which we denote by $M_\Z$, and gives a set of $\binom{rm}{n}-1$ generators
for the $\Z$-module $\bP_n(M_\Z)$.
These will certainly generate $\bP_n(M)$ as an $A$-module, but will be a highly redundant set of generators,
roughly in excess of a factor $r^n$ when $m$ is large (differently from the case $A=\Z$, see the previous subsection).

\bibliography{References}

\def\cprime{$'$} \def\cprime{$'$}
\providecommand{\bysame}{\leavevmode\hbox to3em{\hrulefill}\thinspace}
\providecommand{\MR}{\relax\ifhmode\unskip\space\fi MR }
\providecommand{\MRhref}[2]{%
  \href{http://www.ams.org/mathscinet-getitem?mr=#1}{#2}
}
\providecommand{\href}[2]{#2}
\begin{thebibliography}{GAP07}

\bibitem[Bou89]{Bourbaki:Algebra_I}
Nicolas Bourbaki, \emph{Algebra. {I}. {C}hapters 1--3}, Elements of Mathematics
  (Berlin), Springer-Verlag, Berlin, 1989, Translated from the French, Reprint
  of the 1974 edition. \MR{MR979982 (90d:00002)}

\bibitem[Bou90]{Bourbaki:Algebra_II}
N.~Bourbaki, \emph{Algebra. {II}. {C}hapters 4--7}, Elements of Mathematics
  (Berlin), Springer-Verlag, Berlin, 1990, Translated from the French by P. M.
  Cohn and J. Howie. \MR{MR1080964 (91h:00003)}

\bibitem[CdG07]{CicGra:I}
Serena Cical{\`o} and Willem de~Graaf, \emph{Non-associative {G}r\"obner bases,
  finitely-presented {L}ie rings and the {E}ngel condition}, I{SSAC} 2007, ACM,
  New York, 2007, pp.~100--107. \MR{MR2396190}

\bibitem[CdG09]{CicGra:II}
\bysame, \emph{Non-associative {G}r\"obner bases, finitely-presented {L}ie
  rings and the {E}ngel condition, {II}}, J. Symbolic Comput. \textbf{44}
  (2009), 786--800.

\bibitem[GAP07]{GAP}
The GAP Group, \emph{\textsf{GAP} --- {G}roups, {A}lgorithms, and
  {P}rogramming, version 4.4.10}, 2007, \texttt{http://www.gap-system.org}.

\bibitem[Jac74]{Jac:BAI}
Nathan Jacobson, \emph{Basic algebra. {I}}, W. H. Freeman and Co., San
  Francisco, Calif., 1974. \MR{MR0356989 (50 \#9457)}

\bibitem[Kos90]{Kostrikin:Burnside}
A.~I. Kostrikin, \emph{Around {B}urnside}, Ergebnisse der Mathematik und ihrer
  Grenzgebiete (3) [Results in Mathematics and Related Areas (3)], vol.~20,
  Springer-Verlag, Berlin, 1990, Translated from the Russian and with a preface
  by James Wiegold. \MR{MR1075416 (91i:20038)}

\bibitem[Nic99]{Nickel:Engel}
Werner Nickel, \emph{Computation of nilpotent {E}ngel groups}, J. Austral.
  Math. Soc. Ser. A \textbf{67} (1999), no.~2, 214--222, Group theory.
  \MR{MR1717415 (2000h:20058)}

\bibitem[Zel90]{Zel:Burnside_odd}
E.~I. Zel{\cprime}manov, \emph{Solution of the restricted {B}urnside problem
  for groups of odd exponent}, Izv. Akad. Nauk SSSR Ser. Mat. \textbf{54}
  (1990), no.~1, 42--59, 221. \MR{MR1044047 (91i:20037)}

\end{thebibliography}

\end{document}